\documentclass[10pt,twoside,english]{article}
\usepackage[T1]{fontenc}
\usepackage[latin9]{inputenc}
\usepackage[a4paper]{geometry}
\usepackage{calc}
\usepackage{amsthm}
\usepackage{amsmath}
\usepackage{amssymb}
\usepackage{graphicx}

\makeatletter
\numberwithin{equation}{section}
\numberwithin{figure}{section}

\makeatother

\usepackage{babel}

\def\multiset#1#2{\ensuremath{\left(\kern-.23em\left(\genfrac{}{}{0pt}{}{#1}{#2}\right)\kern-.23em\right)}}

\begin{document}


\title{Extraction of the $i^{th}$ Elementary Symmetric Polynomial \\
from a Product in Binomial Form ${{m_1+\ldots+m_n}\choose i}$
} 

\author{Felix de la Fuente \footnote{mail@felixdelafuente.es}
}

\date{September 17, 2014}

\maketitle

\begin{abstract}
The $i^{th}$ elementary symmetric polynomial of the set of $n$ variables
$\mathcal{R}=\{m_1,m_2,m_3,\ldots,m_n\}$ is isolated from
the expansion of the $i^{th}$ binomial product ${{m_1+\ldots+m_n}\choose i}$ via an alternating sum. 

\end{abstract}


\vphantom{}
\vphantom{}
\vphantom{}
\vphantom{}
\vphantom{}
\vphantom{}

\section{Introduction}
\label{intro}

For a finite root set $\mathcal{R}=\{m_{1},m_{2},m_{3},\ldots,m_{n}\}$
of any positive numbers $m_{i}\in\mathbb{N}$, whether equal or unequal,
the \textbf{$i^{th}$ elementary symmetric polynomial} $e_{i}(n)$ is the sum
of all $i$-tuples of unique elements from $\mathcal{R}$

\begin{eqnarray*}
e_{i}(n) & = & \sum_{1\leqslant j_{1}<j_{2}<\ldots<j_{i}\leqslant n}m_{j_{1}}m_{j_{2}}\ldots m_{j_{i-1}}m_{j_{i}}\\
 & = & \sum_{J\in\binom{1,2,..n}{i}}\prod_{j\in J}m_{j}
\end{eqnarray*}

When the size of the tuples equals $0$, the polynomial accounts for
the empty product $e_{0}(n)=1$.

\vphantom{}

In this way, the set of four numbers $\mathcal{R}=\{m_{1},m_{2},m_{3},m_{4}\}$
generates the following polynomials

\vphantom{}

$e_{0}(4)=1$

$e_{1}(4)=m_{1}+m_{2}+m_{3}+m_{4}$

$e_{2}(4)=m_{1}m_{2}+m_{1}m_{3}+m_{1}m_{4}+m_{2}m_{3}+m_{2}m_{4}+m_{3}m_{4}$

$e_{3}(4)=m_{1}m_{2}m_{3}+m_{1}m_{2}m_{4}+m_{1}m_{3}m_{4}+m_{2}m_{3}m_{4}$

$e_{4}(4)=m_{1}m_{2}m_{3}m_{4}$

\vphantom{}

\begin{figure*}[!hbt]
  \includegraphics[width=1\textwidth]{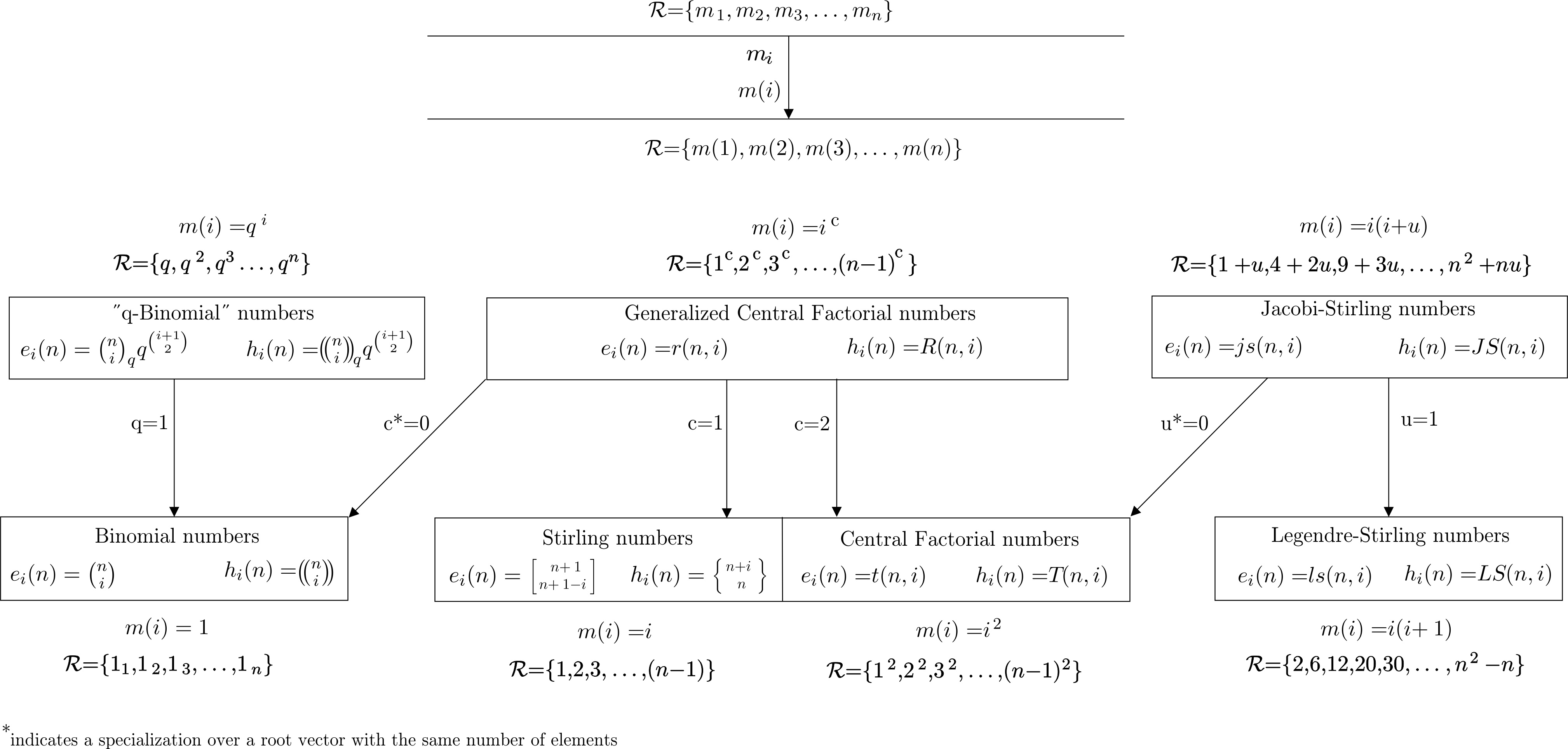}

\caption{Some specialized root sets generating orthogonal families of numbers as elementary and homogeneous symmetric polynomials.}
\label{fig:1}       
\end{figure*}

Algebraically, elementary symmetric polynomials constitute a basis
of the space of symmetric polynomials, and so do others functions
such as homogeneous symmetric polynomials and power sum symmetric
polynomials, which will not be discussed here. Relations among bases
is a classical subject dating back to Guirard, Waring, Newton and Euler \cite{RefA} (see \cite{RefB} for a general overview) as it is the specialization due to the insertion of particular numbers, or $q$-numbers, within the root set. Together
with the homogeneous symmetric polynomials, many of these specializations
have historically been objects of interest in number theory and combinatorics for diverse authors,
and still are focus of current research. Some noteworthy examples are illustrated in Figure~\ref{fig:1}.

\vphantom{}

Let ${N \choose i}=\frac{N(N-1)\cdots(N-i+1)}{i!}$ and $\multiset{N}{i}=\frac{N(N+1)\cdots(N+i-1)}{i!}$
be \textbf{binomial numbers }of the\textbf{ first }and the \textbf{second
kind} respectively.

\vphantom{}

If we consider $N$ as the accumulate of the root set
$N=\sum_{j=1}^{n}m_{j}$, the former binomial numbers represent the \textbf{$i^{th}$ binomial product} $\tbinom{m_{1}+\ldots+m_{n}}{i}=\frac{(m_{1}+\ldots+m_{n})(m_{1}+\ldots+m_{n}-1)\cdots(m_{1}+\ldots+m_{n}-i+1)}{i!}$
for $n>i$.

\vphantom{}

This article presents an expression for the elementary symmetric polynomials of a root set $\mathcal{R}=\{m_1,m_2,m_3,\ldots,m_n\}$ as a function of its associated binomial products. Namely

\begin{eqnarray*}
e_{1}(n)=\tbinom{m_{1}+\ldots+m_{n}}{1}\\
e_{2}(n)=\tbinom{m_{1}+\ldots+m_{n}}{2} & -\multiset{n-1}{0} & \left[\tbinom{m_{1}}{2}+\tbinom{m_{2}}{2}+\tbinom{m_{3}}{2}+\tbinom{m_{4}}{2}+\cdots+\tbinom{m_{n}}{2}\right]\\
e_{3}(n)=\tbinom{m_{1}+\ldots+m_{n}}{3} & -\multiset{n-2}{0} & \left[\tbinom{m_{1}+m_{2}}{3}+\tbinom{m_{1}+m_{3}}{3}+\tbinom{m_{2}+m_{3}}{3}+\cdots\right]\\
 & +\multiset{n-2}{1} & \left[\tbinom{m_{1}}{3}+\tbinom{m_{2}}{3}+\tbinom{m_{3}}{3}+\tbinom{m_{4}}{3}+\cdots+\tbinom{m_{n}}{3}\right]\\
e_{4}(n)=\tbinom{m_{1}+\ldots+m_{n}}{4} & -\multiset{n-3}{0} & \left[\tbinom{m_{1}+m_{2}+m_{3}}{4}+\tbinom{m_{1}+m_{2}+m_{4}}{4}+\tbinom{m_{1}+m_{3}+m_{4}}{4}+\cdots\right]\\
 & +\multiset{n-3}{1} & \left[\tbinom{m_{1}+m_{2}}{4}+\tbinom{m_{1}+m_{3}}{4}+\tbinom{m_{2}+m_{3}}{4}+\cdots\right]\\
 & -\multiset{n-3}{2} & \left[\tbinom{m_{1}}{4}+\tbinom{m_{2}}{4}+\tbinom{m_{3}}{4}+\tbinom{m_{4}}{4}+\cdots+\tbinom{m_{n}}{4}\right]\\
e_{5}(n)=\tbinom{m_{1}+\ldots+m_{n}}{5} & -\multiset{n-4}{0} & \left[\tbinom{m_{1}+m_{2}+m_{3}+m_{4}}{5}+\tbinom{m_{1}+m_{2}+m_{3}+m_{5}}{5}+\tbinom{m_{1}+m_{3}+m_{4}+m_{5}}{5}+\cdots\right]\\
 & +\multiset{n-4}{1} & \left[\tbinom{m_{1}+m_{2}+m_{3}}{5}+\tbinom{m_{1}+m_{2}+m_{4}}{5}+\tbinom{m_{1}+m_{3}+m_{4}}{5}+\cdots\right]\\
 & -\multiset{n-4}{2} & \left[\tbinom{m_{1}+m_{2}}{5}+\tbinom{m_{1}+m_{3}}{5}+\tbinom{m_{2}+m_{3}}{5}+\cdots\right]\\
 & +\multiset{n-4}{3} & \left[\tbinom{m_{1}}{5}+\tbinom{m_{2}}{5}+\tbinom{m_{3}}{5}+\tbinom{m_{4}}{5}+\cdots+\tbinom{m_{n}}{5}\right]
\end{eqnarray*}

\newpage


\section{Extraction of the $i^{th}$ Elementary Symmetric Polynomial}
\label{sec:1}

\vphantom{}

\fbox{\begin{minipage}[t]{1\columnwidth}%
Theorem. \emph{Let $\mathcal{R}=\{m_{1},m_{2},m_{3},\ldots,m_{n}\}$ be a root
set of any positive numbers $m_{j}\in\mathbb{N}$. Then, the $i^{th}$elementary
symmetric polynomial of size $n$ is given by}

\[
e_{i}(n)=\binom{\sum_{k=1}^{n}m_{k}}{i}-\sum_{h=1}^{i-1}(-1)^{h-1}\multiset{n-i+1}{h-1}\left[\sum_{_{J\in\binom{1,2,..n}{i-h}}}\binom{\sum_{_{j\in J}}m_{j}}{i}\right]
\]

\vphantom{}%
\end{minipage}}

\vphantom{}
\vphantom{}

\begin{proof}

Consider the root set of $n$ elements $\mathcal{R}=\{m_{1},m_{2},m_{3},\ldots,m_{n}\}$
and a subset of $s$ elements $\mathcal{Q}_{s}=\{m_{j_{1}},m_{j_{2}},\ldots,m_{j_{s}}\},\,\mathcal{Q}_{s}\in\mathcal{R}$.
The expansion of the product in binomial form $\tbinom{m_{1}+\ldots+m_{n}}{i}=\frac{(m_{1}+\ldots+m_{n})(m_{1}+\ldots+m_{n}-1)\cdots(m_{1}+\ldots+m_{n}-i+1)}{i!}$
can be organized by association to sets of sums involving products
whose factors belong to certain subset $\mathcal{Q}_{s}$.

\vphantom{}

\emph{Example:} Consider the expansion of the $4^{th}$
binomial product $\tbinom{m_{1}+\ldots+m_{n}}{4}$. The summands involving
products of elements belonging to a general subset of size two $Q_{2}=\{m_{j_{1}},m_{j_{2}}\}$
are

\emph{\hspace{1cm}}$+\left[-\frac{22}{4!}m_{j_{1}}m_{j_{2}}\right]$

\emph{\hspace{1cm}}$+\left[\frac{18}{4!}\left(m_{j_{1}}\right)^{2}m_{j_{2}}+\frac{18}{4!}m_{j_{1}}\left(m_{j_{2}}\right)^{2}\right]$

\emph{\hspace{1cm}}$+\left[-\frac{4}{4!}\left(m_{j_{1}}\right)^{3}m_{j_{2}}-\frac{6}{4!}\left(m_{j_{1}}\right)^{2}\left(m_{j_{2}}\right)^{2}-\frac{4}{4!}m_{j_{1}}\left(m_{j_{2}}\right)^{3}\right]$

\vphantom{}

The sum of the exponents of the factors of every product add up to
an accumulate power $p$ which ranges from the size of the set $s$
to the maximal value $i$. The repetitions of every product are given
by the coefficient $K_{p}^{\lambda(p)_{k}}$, where $\lambda(p)_{k}$is
the ordered partition of $p$ into $i$ parts corresponding to the
distribution of the exponents. Extending the former example for the
expansion of the general $i^{th}$ binomial product $\tbinom{m_{1}+\ldots+m_{n}}{i}$,
the summands involving factors from $\mathcal{Q}_{2}$ are

\vphantom{}

\emph{\hspace{1cm}}$+\left[K_{2}^{(1,1)}m_{j_{1}}m_{j_{2}}\right]$

\emph{\hspace{1cm}}$+\left[K_{3}^{(2,1)}\left(m_{j_{1}}\right)^{2}m_{j_{2}}+K_{3}^{(1,2)}m_{j_{1}}\left(m_{j_{2}}\right)^{2}\right]$

\emph{\hspace{1cm}}$+\left[K_{4}^{(3,1)}\left(m_{j_{1}}\right)^{3}m_{j_{2}}+K_{4}^{(2,2)}\left(m_{j_{1}}\right)^{2}\left(m_{j_{2}}\right)^{2}+K_{4}^{(1,3)}m_{j_{1}}\left(m_{j_{2}}\right)^{3}\right]$

\emph{\hspace{1cm}}$\,\vdots$

\emph{\hspace{1cm}}$+\left[\sum_{k}K_{p}^{\lambda(p)_{k}}(m_{j_{1}}m_{j_{2}})_{p}^{\lambda(p)_{k}}\right]$

\emph{\hspace{1cm}}$\,\vdots$

\emph{\hspace{1cm}}$+\left[\sum_{k}K_{i}^{\lambda(i)_{k}}(m_{j_{1}}m_{j_{2}})_{i}^{\lambda(i)_{k}}\right]$

\vphantom{}

\emph{\hspace{1cm}}$=\sum_{p=2}^{i}\left[\sum_{k}K_{p}^{\lambda(p)_{k}}(m_{j_{1}}m_{j_{2}})_{p}^{\lambda(p)_{k}}\right]$

\vphantom{}

\vphantom{}

We will make now two key observations, the first one is that for any
subset $\mathcal{Q}_{s}=\{m_{j_{1}},m_{j_{2}},\ldots,m_{j_{s}}\}\in\mathcal{R}=\{m_{1},m_{2},m_{3},\ldots,m_{n}\}$
and a fixed value of $i$, the $i^{th}$ binomial expansion over $\mathcal{R}$,
$\tbinom{m_{1}+\ldots+m_{n}}{i}$, will generate \emph{the same} coefficient
$K_{p}^{\lambda(p)_{k}}$ for the product of elements of $\mathcal{Q}_{s}$
corresponding to such partition, and it also does over any subset
of intermediate size $s<t<n$, $\mathcal{Q}_{t}=\{m_{j_{1}},m_{j_{2}},\ldots,m_{j_{t}}\}$.

\vphantom{}

This allows a convenient simplification of the nomenclature since
these coefficients will not affect the forthcoming sieve: in the former
example $\mathcal{Q}_{2}=\{m_{j_{1}},m_{j_{2}}\}$, the sum of all
products of two distinct elements arosen to all possible powers, that
is the sums of products over subsets of size two, will be referred
more compactly as $\left[2\, tuples\right]^{2\leqslant p\leqslant i}$
of just $2\, tuples$.

\vphantom{}

$\sum_{\overset{j\in J}{J\in\binom{1,...,n}{2}}}\sum_{p=2}^{i}\left[\sum_{k}K_{p}^{\lambda(p)_{k}}(m_{j_{1}}m_{j_{2}})_{p}^{\lambda(p)_{k}}\right]\longmapsto\left[2-tuples\right]^{2\leqslant p\leqslant i}\mbox{ or }2\, tuples$ 

\vphantom{}

On the other hand, it can be observed that the sums corresponding
to the subset of maximal size, that is of the size of the binomial
product $s=i$, have a coefficient $K_{i}^{(\overbrace{1,1,...,1}^{n})}=1$
and the $\left[i\, tuples\right]^{p=i}=\sum_{J\in\binom{1,2,..n}{i}}\prod_{j\in J}m_{j}=e_{i}(n)$
are exactly the \emph{elementary symmetric polynomial} of order $i$.

\vphantom{}

These two facts allow naturally the following construction where coinciding
typologies of $\left[t\, tuples\right]^{t\leqslant p\leqslant i}$,
corresponding to $\mathcal{Q}_{t}\in R$, are ordered by reverse inclusion
horizontally, and exploit afterwards this coincidence to extract $e_{i}(n)$
from $\tbinom{m_{1}+\ldots+m_{n}}{i}$.

\vphantom{}

  \includegraphics[width=.9\textwidth]{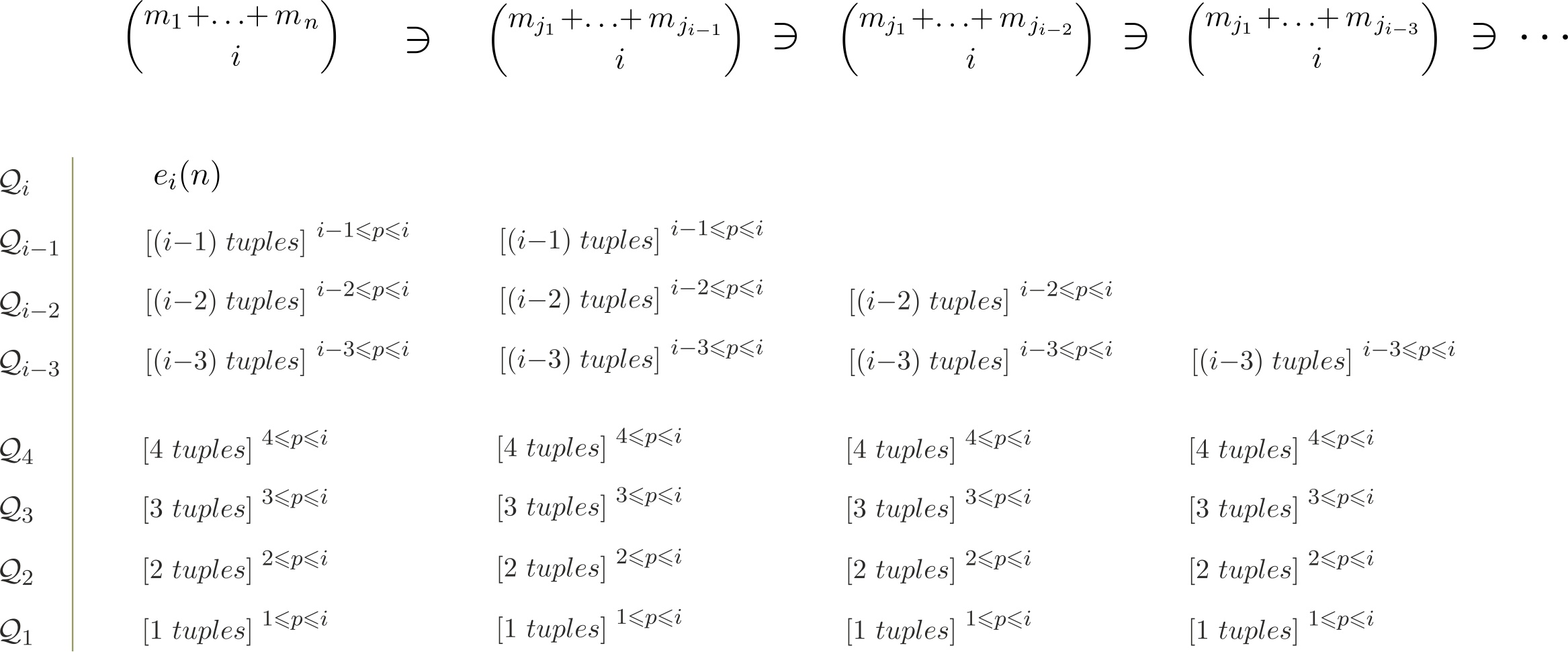}

\vphantom{}

Now, in order to generate all the varieties of $s\, tuples$ of the
first column in any of the others, call it $\binom{m_{j_{1}}+\ldots+m_{j_{s}}}{i}$,
it is necessary to add up all binomials whose numerators cover all
possible sums of size $s$ of elements within the root set $\mathcal{R}=\{m_{1},m_{2},m_{3},\ldots,m_{n}\}$,
that is

\[
\sum_{_{J\in\tbinom{1,\ldots,n}{s}}^{j\in J}}\binom{\sum m_{j}}{i}
\]

Every column has exactly the same typologies of $t\, tuples$,
but not the same quantity except for $t=s$, which appear once. Indeed,
observe that all $t\, tuples$ for $1\leqslant t<s$ will be repeated
as many times as sets of size $p-t$ can be formed with the $n$ elements
of the root set, from which the $t$ that remain fixed are excluded,
and in consequence every set of $t\, tuples$ appear with multiplicity
$\binom{n-t}{s-t}$ .


  \includegraphics[width=.5\textwidth]{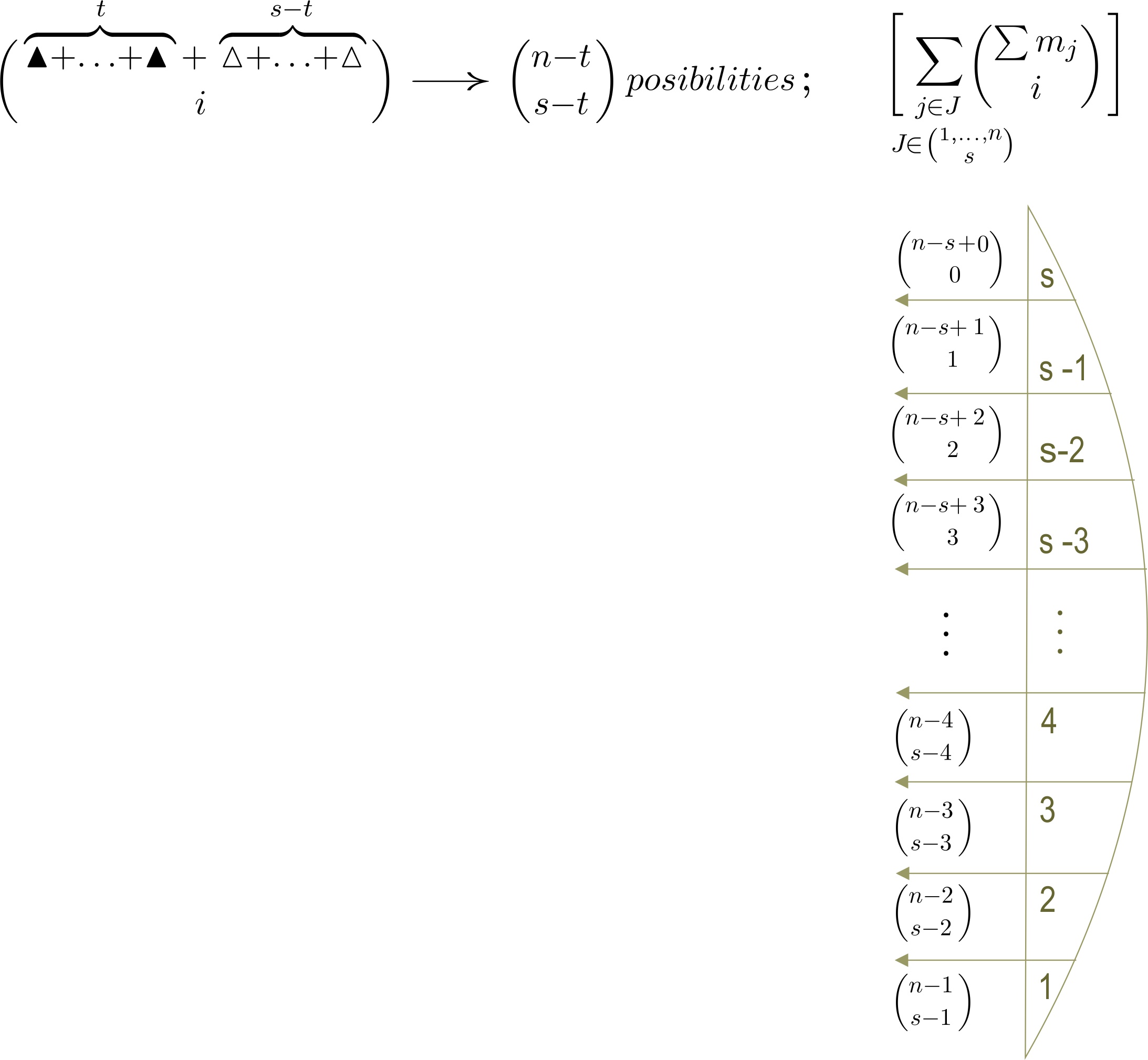}

\vphantom{}
\newpage

Finally, all sums contain the same typologies of elements that the
others and the same typologies that the expansion of $\tbinom{m_{1}+\ldots+m_{n}}{i}$
whose multiplicity is $1$ for all $t\, tuples$; the elementary symmetric
polynomial $e_{i}(n)$ can be then isolated by extracting all $t\, tuples$
for $1\leqslant t\leqslant i-1$ of the other columns weighted by
the adequate coefficients $C_{h}$

\[
e_{i}(n)=\binom{\sum_{i=1}^{n}m_{i}}{i}-\left[\sum_{h=1}^{i-1}C_{h}\left[\sum_{_{J\in\binom{1,2,..n}{i-h}}}\binom{\sum_{_{j\in J}}m_{j}}{i}\right]\right]
\]

\vphantom{}

  \includegraphics[width=1\textwidth]{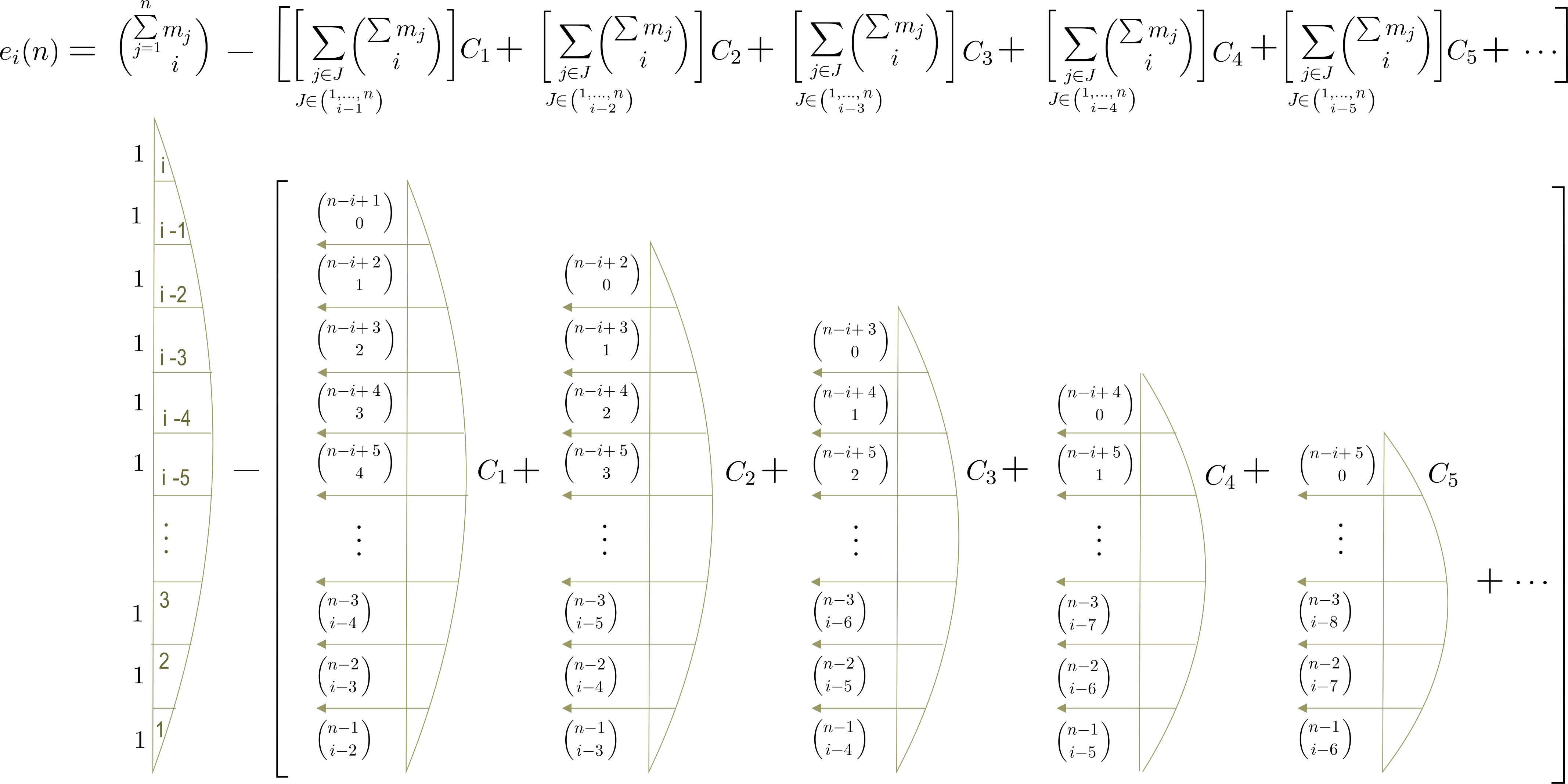}

\vphantom{}

Since the typologies of the objects are the same, suffice it to consider
the multiplicities to tackle the computation of the $C_{h}$ and hence
the following series of expressions can be constructed:

\vphantom{}

$1=\tbinom{n-i+1}{0}C_{1}$

$1=\tbinom{n-i+2}{1}C_{1}+\tbinom{n-i+2}{0}C_{2}$

$1=\tbinom{n-i+3}{2}C_{1}+\tbinom{n-i+3}{1}C_{2}+\tbinom{n-i+3}{0}C_{3}$

$1=\tbinom{n-i+4}{3}C_{1}+\tbinom{n-i+4}{2}C_{2}+\tbinom{n-i+4}{1}C_{3}+\tbinom{n-i+4}{0}C_{4}$

$\,\vdots$

$1=\tbinom{n-i+h}{h-1}C_{1}+\tbinom{n-i+h}{h-2}C_{2}+\tbinom{n-i+h}{h-3}C_{3}+\tbinom{n-i+h}{h-4}C_{4}+\cdots+\tbinom{n-i+h}{1}C_{h-1}+\tbinom{n-i+h}{0}C_{h}$

\vphantom{}

where the $h^{th}$ term can be defined as the \textbf{recurrence}
of degree $h-1$:

\[
C_{h}=1-\sum_{k=1}^{h-1}C_{k}\binom{n-i+h}{h-k}
\]

or otherwise be inserted in the \textbf{complete convolution}

\begin{equation}\label{eq:1}
1=\sum_{k=1}^{h}C_{k}\binom{n-i+h}{h-k}
\end{equation}

\vphantom{}

We sill resolve the sequence of coefficients $C_{1},C_{2},C_{3},\ldots,C_{h}$
twice: first in a very immediate but scarcely elucidating way by considering
a degeneration of the Vandermonde's identity and then, alternatively, in
analytical terms through formal power series and generating functions.

\begin{enumerate}
\item Consider the complete convolution (\ref{eq:1}) for the term $C_{h+1}$\\
\[
1=\sum_{k=1}^{h+1}C_{k}\binom{n-i+h+1}{h+1-k}=\sum_{k=0}^{h}C_{k+1}\binom{n-i+h+1}{h-k}
\]
\\
and the Vandermonde's identity $\binom{l+m}{h}=\sum_{k=0}^{h}\binom{l}{k}\binom{m}{h-k}$.
The first expression can be regarded as a ``limit'' specialization
of the second one for the values $l=-n+i-1$ and $m=n-i+h+1$
\[
\binom{h}{h}=1=\sum_{k=0}^{h}\binom{-n+i-1}{k}\binom{n-i+h+1}{h-k}
\]
\\
From where the coefficients can be pairwise identified with the numbers\\
\begin{eqnarray*}
C_{k} & = & \binom{-n+i-1}{k-1}\\
 & = & (-1)^{k-1}\multiset{n-i+1}{k-1}
\end{eqnarray*}

\item Let us consider again the complete convolution (\ref{eq:1}) for the term $C_{h}$ and associate
formal power series to the terms of the sequence
\begin{eqnarray*}
1 & = & \sum_{k=1}^{h}C_{k}\binom{n-i+h}{h-k}\\
\sum_{h=1}^{\infty}x^{h-1} & = & \sum_{h=1}^{\infty}\left[\sum_{k=1}^{h}C_{k}\binom{n-i+h}{h-k}\right]x^{h-1}\\
 & = & \sum_{k=1}^{h}C_{k}\sum_{h=k}^{\infty}\binom{n-i+h}{h-k}x^{h-1}
\end{eqnarray*}
\\
with the infinite geometric series on the left hand side and the binomial
series on the right hand side, 
\begin{eqnarray*}
\frac{1}{1-x} & = & \sum_{k=1}^{h}C_{k}\frac{x^{k-1}}{(1-x)^{n-i+k+1}}\\
 & = & \frac{1}{x(1-x)^{n-i+1}}\sum_{k=1}^{h}C_{k}\left(\frac{x}{(1-x)}\right)^{k}\\
 & \mbox{and }\\
x(1-x)^{n-i} & = & \sum_{k=1}^{h}C_{k}\left(\frac{x}{(1-x)}\right)^{k}
\end{eqnarray*}
\\
performing the transformation $x\longmapsto\frac{x}{(1+x)}$ the equality
yields a form that admits a term by term identification with the expansion
of the left hand side function
\begin{eqnarray*}
\frac{x}{(1+x)^{n-i+1}} & = & \sum_{k=1}^{h}C_{k}x^{k}\\
\sum_{k=0}^{\infty}(-1)^{k}\binom{n-i+k+1}{k}x^{k+1} & = & \sum_{k=0}^{\infty}C_{k+1}x^{k+1}
\end{eqnarray*}
\\
from where 
\begin{eqnarray*}
C_{k} & = & (-1)^{k-1}\binom{n-i+k}{k-1}\\
 & = & \binom{-n+i-1}{k-1}\\
 & = & (-1)^{k-1}\multiset{n-i+1}{k-1}
\end{eqnarray*}

\end{enumerate}

\end{proof}



\end{document}